\documentclass[a4paper,12pt]{article}
\usepackage[cp1251]{inputenc}
\usepackage[english]{babel}
\usepackage[tbtags]{amsmath}
\usepackage{amsfonts,amssymb,mathrsfs,amscd}
\usepackage{wrapfig}
\usepackage[dvips]{graphicx}
\newtheorem{theorem}{Theorem}[section]
\newtheorem{lemma}{Lemma}[section]
\newtheorem{definition}{Definition}[section]

\date{}
 \numberwithin{equation}{section}

\begin{document}
 \large
 \centerline{\bf N-valid trees in wavelet theory on  Vilenkin groups \footnote[1]{This research was carried out with the financial
support of  the Russian Foundation for Basic Research (grant
no.~13-01-00102).} }
 \centerline{\bf G.\,S.Berdnikov, S.\,F.~Lukomskii}

\noindent
 N.G.\ Chernyshevskii Saratov State University\\
  evrointelligent@gmail.com\\ LukomskiiSF@info.sgu.ru\\
 MSC:Primary 42C40; Secondary 11R56, 43A70\\


\begin{abstract}
 We consider a class of
 $(N,M)$-elementary step functions on the $p$-adic Vilenkin group.
 We prove that $(N,M)$-elementary step
 function generates a MRA on $p$-adic Vilenkin group iff it is
 generated by a special $N$-valid  rooted tree on the set of vertices
 $\{0,1,\dots p-1\}$ with the vector $(0,...,0)\in \mathbb Z^N$ as a root.
Bibliography: 15 titles.
\end{abstract}
\noindent
 keywords: zero-dimensional group, Vilenkin group, multiresolution
analysis, wavelet bases, tree.


\section{Introduction}\label{s1}

  In articles \cite{1}-\cite{4} first
 examples of orthogonal wavelets on the dyadic Cantor  group ($p=2$) are
 constructed and their properties are studied.
  Yu.Farkov \cite{5}-\cite{7} found necessary and
 sufficient conditions for a refinable function to generate an orthogonal
 MRA in the $L_2(\mathfrak G)$
 -spaces on the $p$-adic Vilenkin group $\mathfrak G$. These conditions use the
 Strang-Fix and the modified Cohen  properties.

In \cite{7}  this
 construction  is  given in a concrete fashion for p = 3.
 In \cite{8}, some algorithms for constructing orthogonal and biorthogonal compactly
 supported wavelets on Vilenkin groups are proposed. In \cite{5}-\cite{8} two
 types of orthogonal wavelet examples are constructed: step
 functions and sums of Vilenkin  series.

 Khrennikov, Shelkovich, and Skopina
\cite{10},\cite{11} introduced the concept of a~$p$-adic MRA with
orthogonal refinable function, and described a~general pattern for
their construction. This method was developed for an orthogonal
refinable function $\varphi$ with condition ${\rm
supp}\widehat{\varphi}\subset B_0(0)$, where $B_0(0) = \{x :
|x|_p\le 1 \}$ is the unit ball in the field $\mathbb Q_p$.
Similar results were obtained for an arbitrary zero-dimensional
group \cite{13}. The condition  ${\rm
supp}\widehat{\varphi}\subset B_0(0)$ is very important. S.
Albeverio, S. Evdokimov, M. Skopina \cite{12} proved that if a
refinable step function $\varphi$ generates an orthogonal $p$-adic
MRA, then ${\rm supp}\hat \varphi(\chi)\subset B_0(0)$.

On the other hand on Vilenkin groups Yu.A.Farkov constructs
examples of step refinable functions $\varphi$, which generate an
orthogonal MRA with ${\rm supp}\widehat{\varphi}\subset G_1^\bot$.
In the author's work \cite{14} a necessary condition for a support
of orthogonal refinable step function are found: if step refinable
 $(1,M)$-elementary function $\varphi$ generates an orthogonal MRA
 on $p$-adic Vilenkin group, then ${\rm supp}\widehat{\varphi}\subset
 G_{p-2}^\bot$. In \cite{15} some
 trees was used to construct refinable function .

 In this work we consider more general situation and
 study a structure of the set
 ${\rm supp}\widehat{\varphi}$.  We define a concept of $N$-valid tree and  prove
 that $(N,M)$-elementary function $\varphi$
 generates an orthogonal MRA
 on $p$-adic Vilenkin group iff the function $\varphi$  is
 generated by means of some $N$-valid tree. For any $N$-valid tree we give an
 algorithm for constructing corresponding refinable
 function and orthogonal wavelets.

The paper is organized as follows. We  consider $p$-adic Vilenkin
group $\mathfrak G$ as a zero-dimensional group $(G,\dot + )$ with
condition  $pg_n$=0. Therefore, in section 2, we recall some
concepts and facts from the theory of  zero-dimensional group. We
will systematically use the notation and the results from
\cite{13},\cite{14}.

In section 3 and the following sections we consider MRA on
$p$-adic Vilenkin group $\mathfrak G$. In section 3 we study
refinable step-functions which generate the orthogonal MRA. We
define a class of $(N,M)$-elementary set and prove that the shifts
system $\varphi (x\dot - h)_{h\in H_0}$ is orthonormal if ${\rm
supp}\widehat{\varphi}$ is $(N,M)$-elementary set.

In section 4 we introduce such  concepts as "a set generated by a
tree" and "a refinable step function generated by a tree" and
prove, that any rooted $N$-valid tree generates a refinable step
function that generates an orthogonal MRA on Vilenkin group.

In section 5 we give an algorithm for constructing orthogonal
wavelets according to the tree.

\section{Preliminaries}\label{s2}
We will consider the Vilenkin group as a  locally compact
zero-dimensional Abelian group
 with additional condition $p_ng_n=0$. Therefore we start with some basic notions
  and facts related to analysis on zero-dimensional groups. One may find more information on the topic in \cite{12}--\cite{14}.

Let $(G,\dot + )$~be a~locally compact zero-dimensional Abelian
group with the topology generated by a~countable system of open
subgroups
$$
\cdots\supset G_{-n}\supset\cdots\supset G_{-1}\supset G_0\supset
G_1\supset\cdots\supset G_n\supset\cdots
$$
where
$$
\bigcup_{n=-\infty}^{+\infty}G_n= G,\quad
 \quad \bigcap_{n=-\infty}^{+\infty}G_n=\{0\},
$$
 $p_n$ is an order of quotient group $G_n/G_{n+1}$.
 We will always
assume that all~$p_n$ are prime numbers. We will name such chain
as \it basic chain. \rm  In this case, a~base of the topology is
formed by all possible cosets~$G_n\dot + g$, $g\in G$.

We further define the numbers $(\mathfrak
m_n)_{n=-\infty}^{+\infty}$ as follows:
$$
\mathfrak m_0=1,\qquad \mathfrak m_{n+1}=\mathfrak m_n\cdot p_n.
$$
Let~$\mu$ be a Haar measure on~$G$, we know that $\mu
G_n=\frac{1}{\mathfrak m_n}$. Further, let
$\smash[b]{\displaystyle\int_{G}f(x)\,d\mu(x)}$~be the absolutely
convergent integral of the measure~$\mu$.

Given $n\in\mathbb Z$, consider an element $g_n\in G_n\setminus
G_{n+1}$ and fix~it. Then any $x\in G$ has a~unique representation
in the form
\begin{equation}
\label{eq1.1} x=\sum_{n=-\infty}^{+\infty}a_ng_n, \qquad
a_n=\overline{0,p_n-1}.
\end{equation}
The sum~\eqref{eq1.1} contain finite number of terms with negative
subscripts, that~is,
\begin{equation}
\label{eq1.2} x=\sum_{n=m}^{+\infty}a_ng_n, \qquad
a_n=\overline{0,p_n-1}, \quad a_m\ne 0.
\end{equation}
We will name system $(g_n)_{n\in \mathbb Z}$ as {\it a basic
system}.

Classical examples of zero-dimensional groups are Vilenkin groups
and groups of $p$-adic numbers (see~\cite[Ch.~1, \S\,2]{12}).
 A direct sum of cyclic groups $Z(p_k)$ of order~$p_k$,
 $k\in\mathbb Z$, is called a~\textit{Vilenkin group}. This means
 that the elements of a~Vilenkin group are infinite sequences
 $x=(x_k)_{k=-\infty}^{+\infty}$ such that:
 \begin{itemize}
 \item[1)] $x_k=\overline{0,p_k-1}$; \item[2)] only a~finite number
 of~$x_k$ with negative subscripts are different from zero;
 \item[3)] the group operation~$\dot+ $ is the coordinate-wise
 addition modulo~$p_k$, that~is,
 $$
 x\dot+ y=(x_k\dot+ y_k), \qquad x_k\dot+ y_k=(x_k+y_k)\ \
 \operatorname{mod}p_k.
 $$
 \end{itemize}
 A topology on such group is generated by the chain of subgroups
 $$
  G_n=\bigl\{x\in G:x=(\dots,0,0,\dots,0,x_n,x_{n+1},\dots),\
   x_\nu=\overline{0,p_\nu-1},\ \nu\ge n\bigr\}.
 $$
 The elements $g_n=(\dots,0,0,1,0,0,\dots)$ form a basic system.
 From definition of the operation $\dot+$ we have $p_ng_n=0$.
 Therefore we will name a zero-dimensional group $(G,\dot+)$ with the condition $p_ng_n=0$ as Vilenkin group.

By $X$ we denote  the collection
 of the characters of a~group $(G,\dot+ )$; it is
a~group with respect to multiplication, too. Also let
$G_n^\bot=\{\chi\in X:\forall\,x\in G_n\  , \chi(x)=1\}$ be the
annihilator of the group~$G_n$. Each annihilator~$ G_n^\bot$ is
a~group with respect to multiplication, and the subgroups~$
G_n^\bot$ form an~increa\-sing sequence
\begin{equation}
\label{eq1.3} \cdots\subset G_{-n}^\bot\subset\cdots\subset
G_0^\bot \subset G_1^\bot\subset\cdots\subset
G_n^\bot\subset\cdots
\end{equation}
with
$$
\bigcup_{n=-\infty}^{+\infty} G_n^\bot=X \quad\text{and} \quad
\bigcap_{n=-\infty}^{+\infty} G_n^\bot=\{1\},
$$
the quotient group $ G_{n+1}^\bot/ G_n^\bot$ having order~$p_n$.
The group of characters~$X$ is a zero-dimensional group with a
basic chain \eqref{eq1.3}. The group  may be supplied with the
topology using the chain of subgroups~\eqref{eq1.3}, the family of
the cosets $ G_n^\bot\cdot\chi$, $\chi\in X$, being taken
as~a~base of the topology. The collection of such cosets, along
with the empty set, forms the~semiring~${\mathscr X}$. Given
a~coset $ G_n^\bot\cdot\chi$, we define a~measure~$\nu$ on it by
$\nu( G_n^\bot\cdot\chi)=\nu( G_n^\bot)= \mathfrak m_n$ (so that
always $\mu( G_n)\nu( G_n^\bot)=1$). The measure~$\nu$ can be
extended onto the $\sigma$-algebra of measurable sets in the
standard way. One then forms the absolutely convergent integral
$\displaystyle\int_XF(\chi)\,d\nu(\chi)$ using this measure.

The value~$\chi(g)$ of the character~$\chi$ at an element $g\in G$
will be denoted by~$(\chi,g)$. The Fourier transform~$\widehat f$
of an~$f\in L_2( G)$  is defined~as follows
$$
\widehat f(\chi)=\int_{ G}f(x)\overline{(\chi,x)}\,d\mu(x)=
\lim_{n\to+\infty}\int_{ G_{-n}}f(x)\overline{(\chi,x)}\,d\mu(x),
$$
with the limit being in the norm of $L_2(X)$. For any~$f\in L_2(G)$,
the inversion formula is valid
$$
f(x)=\int_X\widehat f(\chi)(\chi,x)\,d\nu(\chi)
=\lim_{n\to+\infty}\int_{ G_n^\bot}\widehat
f(\chi)(\chi,x)\,d\nu(\chi);
$$
here the limit also signifies the convergence in the norm of~$L_2(
G)$. If $f,g\in L_2( G)$ then the Plancherel formula is valid
$$
\int_{ G}f(x)\overline{g(x)}\,d\mu(x)= \int_X\widehat
f(\chi)\overline{\widehat g(\chi)}\,d\nu(\chi).
$$
\goodbreak

Provided with this topology, the group of characters~$X$ is
a~zero-dimensional locally compact group; there is, however,
a~dual situation: every element $x\in G$ is a~character of the
group~$X$, and~$ G_n$ is the annihilator of the group~$ G_n^\bot$.
The union of disjoint sets $E_j$ we will denote by $\bigsqcup
E_j$.

 For any $n\in \mathbb Z$ we choose a character $r_n\in  G_{n+1}^{\bot}\backslash G_n^{\bot}$
 and fixed it. $(r_n)_{n\in \mathbb Z}$ is called a Rademacher system. Let us denote
  $$
    H_0=\{h\in G: h=a_{-1}g_{-1}\dot+a_{-2}g_{-2}\dot+\dots \dot+ a_{-s}g_{-s}, s\in \mathbb
    N,\ a_j=\overline{0,p-1}\},
  $$
  $$
    H_0^{(s)}=\{h\in G: h=a_{-1}g_{-1}\dot+a_{-2}g_{-2}\dot+\dots \dot+
    a_{-s}g_{-s},\ a_j=\overline{0,p-1}
    \},s\in \mathbb N.
  $$
  The set $H_0$ is an analog of the set $\mathbb N_0=\mathbb N\bigsqcup \{0\}$.

 If in the zero-dimensional group $G\ p_n=p$ for any $n\in\mathbb Z$
 then we can define the mapping ${\cal
 A}\colon G\to G$ by
 ${\cal A}x:=\sum_{n=-\infty}^{+\infty}a_ng_{n-1}$, where
 $x=\sum_{n=-\infty}^{+\infty}a_ng_n\in G$.  The mapping~${\cal A}$ is called
 a  dilation operator if~${\cal A}(x\dot+ y)={\cal A}x\dot + {\cal A}y$ for all
 $x,y\in G$. By definition, put $(\chi {\cal A},x)=(\chi, {\cal
 A}x)$.

  \begin{lemma}[\cite{14}]
 For any zero-dimensional group\\
 1) $\int\limits_{G_0^\bot}(\chi,x)\,d\nu(\chi)={\bf 1}_{G_0}(x)$,
 2) $\int\limits_{G_0}(\chi,x)\,d\mu(x)={\bf 1}_{G_0^\bot}(\chi)$.\\
 \end{lemma}

 \begin{lemma}[\cite{14}]
  If $p_n=p$ for any $n\in \mathbb Z$ and the mapping ${\cal A}$ is additive then \\
  1) $\int\limits_{G_n^\bot}(\chi,x)\,d\nu(\chi)=p^n{\bf
  1}_{G_n}(x)$,\\
  2) $\int\limits_{G_n}(\chi,x)\,d\mu(x)=\frac{1}{p^n}{\bf
  1}_{G_n^\bot}(\chi)$.
  \end{lemma}
\begin{lemma}[\cite{14}]
Let  $\chi_{n,s}=r_n^{\alpha_n}r_{n+1}^{\alpha_{n+1}}\dots
r_{n+s}^{\alpha_{n+s}}$ be a character which does not belong to
$G_n^\bot$. Then
$$
\int\limits_{G_n^\bot\chi_{n,s}}(\chi,x)\,d\nu(\chi)=p^n(\chi_{n,s},x){\bf
1}_{G_n}(x).
$$
\end{lemma}
\begin{lemma}[\cite{14}]
Let
$h_{n,s}=a_{n-1}g_{n-1}\dot+a_{n-2}g_{n-2}\dot+\dots\dot+a_{n-s}g_{n-s}\notin
G_n$. Then
$$
\int\limits_{G_n\dot+h_{n,s}}(\chi,x)\,d\mu(x)=\frac{1}{p^n}(\chi,h_{n,s}){\bf
1}_{G_n^\bot}(\chi).
$$
\end{lemma}
 \begin{definition}[\cite{14}]
Let $M,N\in\mathbb N$.
We denote by  ${\mathfrak D}_M(G_{-N})$ the set of step-functions
 $f\in L_2(G)$ such that 1) ${\rm supp}\,f\subset G_{-N}$, and 2)
 $f$ is constant on cosets $G_M\dot+g$. ${\mathfrak
 D}_{-N}(G_{M}^\bot)$ is defined similarly.
 \end{definition}
\begin{lemma}[\cite{14}]
 Let $M,N\in\mathbb N$. $f\in \mathfrak D_M(G_{-N})$ if and only if $\hat f\in \mathfrak
 D_{-N}(G_M^\bot)$.
\end{lemma}
\section{MRA and refinable function on Vilenkin groups}\label{s3}

 In what follows we will consider groups $G$ for which
$p_n=p$
  and $pg_n=0$ for any $n\in \mathbb Z$. We know that it is a
  Vilenkin group. We will denote a Vilenkin group as $\mathfrak G$.

  In this group  we can choose Rademacher functions
  in various ways.
  We define Rademacher functions with the equation
  $$
  \left(r_n,\sum_{k\in\mathbb Z}a_kg_k\right)=\exp\left(\frac{2\pi i}{p}a_n\right).
  $$
  In this case
  $$
  (r_n,g_k)=\exp\left(\frac{2\pi i}{p}\delta_{nk}\right).
  $$
 Our main objective is to
 find a simple  algorithm  to get a refinable step-function that generates an orthogonal MRA on Vilenkin group.
 \begin{definition}
  A family of closed subspaces $V_n$, $n\in\mathbb Z$,
 is said to be a~multi\-resolution analysis  of~$L_2(\mathfrak G)$
 if the following axioms are satisfied:
 \begin{itemize}
 \item[A1)] $V_n\subset V_{n+1}$;
 \item[A2)] ${\vrule width0pt
 depth0pt height11pt} \overline{\bigcup_{n\in\mathbb
 Z}V_n}=L_2(\mathfrak G)$ and $\bigcap_{n\in\mathbb Z}V_n=\{0\}$;
 \item[A3)] $f(x)\in V_n$  $\Longleftrightarrow$ \ $f({\cal A} x)\in V_{n+1}$ (${\cal A}$~is a~dilation
 operator);
 \item[A4)] $f(x)\in V_0$ \ $\Longrightarrow$ \
 $f(x\dot - h)\in V_0$ for all $h\in H_0$; ($H_0$ is analog of $\mathbb
 Z$).
  \item[A5)] there exists
 a~function $\varphi\in L_2(\mathfrak G)$ such that the system
 $(\varphi(x\dot - h))_{h\in H_0}$ is an orthonormal basis
 for~$V_0$.
\end{itemize}

 A function~$\varphi$ occurring in axiom~A5 is called
a~\textit{scaling function}.
\end{definition}

 Next we will follow the conventional approach. Let
 $\varphi(x)\,{\in}\, L_2(\mathfrak G)$, and assume that
 $(\varphi(x\dot -\nobreak h))_{h\in H_0}$ is an~orthonormal
 system in~$L_2(\mathfrak G)$. With the function~$\varphi$ and the
 dilation operator~${\cal A}$, we define the linear subspaces
 $L_n=(\varphi({\cal A}^nx\dot - h))_{h\in H_0}$ and
 closed subspaces $V_n=\overline{L_n}$. It is evident that the functions
  $p^{\frac{n}{2}}\varphi({\cal A^n}x \dot-h)_{h\in H_0}$ form
 an orthonormal basis for $V_n$, $n\in \mathbb Z$.   If subspaces $V_n$ form
 a~MRA, then the function~$\varphi$ is said to \textit{generate}
 an~MRA in~$L_2(\mathfrak G)$. If a function $\varphi$ generates an MRA, then we obtain from the axiom A1
\begin{equation}
  \label{eq3.1}
  \varphi(x)=\sum_{h\in H_0}\beta_h\varphi({\cal
  A}x\dot-h)\;\;\left(\sum|\beta_h|^2<+\infty\right).
 \end{equation}
  Therefore we will look up a~function
 $\varphi\in L_2(\mathfrak G)$, which generates an~MRA
 in~$L_2(\mathfrak G)$, as a~solution of the refinement
 equation (\ref{eq3.1}), A solution of refinement equation (\ref{eq3.1}) is called a {\it refinable function}.
 \begin{lemma}[\cite{14}]
Let $\varphi \in \mathfrak D_M(\mathfrak G_{-N})$ be a solution of
(\ref{eq3.1}). Then
\begin{equation} \label{eq3.2}
\varphi(x)=\sum_{h\in H_0^{(N+1)}}\beta_h\varphi({\cal A}x\dot-h)
 \end{equation}
\end{lemma}
 The refinement equation (\ref{eq3.2}) may be written in the form
 \begin{equation}                                      \label{eq3.3}
 \hat\varphi(\chi)=m_0(\chi)\hat\varphi(\chi{\cal
  A}^{-1}),
 \end{equation}
  where

 \begin{equation}                                      \label{eq3.4}
 m_0(\chi)=\frac{1}{p}\sum_{h\in
 H_0^{(N+1)}}\beta_h\overline{(\chi{\cal A}^{-1},h)}
 \end{equation}
is a mask of the equation (\ref{eq3.3}).
\begin{lemma}[\cite{14}]
Let $\varphi\in\mathfrak D_M(\mathfrak G_{-N})$. Then the mask
$m_0(\chi)$ is constant on cosets $\mathfrak G_{-N}^\bot\zeta$. If
$\hat\varphi(\mathfrak G_{-N}^\bot)\neq 0$ then $m_0(\mathfrak
G_{-N}^\bot)=1$.
 \end{lemma}
  \begin{lemma}[\cite{14}]
 The mask  $m_0(\chi)$ is a periodic  function with any period
 $r_1^{\alpha_1}r_2^{\alpha_2}\dots r_s^{\alpha_s}$ $(s\in\mathbb
 N,\; \alpha_j=\overline{0,p-1},\;j=\overline{1,s})$.
 \end{lemma}
 So, if $m_0(\chi)$ is a mask of (\ref{eq3.3}) then\\
 T1) $m_0(\chi)$
 is constant on cosets  $\mathfrak G_{-N}^\bot\zeta$,\\
 T2) $m_0(\chi)$ is periodic with any period
  $r_1^{\alpha_1}r_2^{\alpha_2}\dots
 r_s^{\alpha_s}$, $\alpha_j=\overline{0,p-1}$, \\
 T3)
 $m_0(\mathfrak G_{-N}^\bot)=1$. \\
 Therefore we will assume that  $m_0$
 satisfies these conditions.

  \begin{theorem}[\cite{14}]
 $m_0(\chi)$ is a mask of  equation (\ref{eq3.3}) on the class $\mathfrak
 D_{-N}(\mathfrak G_M^\bot)$ if and only if
 \begin{equation} \label{eq3.5}
 m_0(\chi)m_0(\chi{\cal A}^{-1})\dots m_0(\chi{\cal A}^{-M-N})=0
 \end{equation}
 on $\mathfrak G_{M+1}^\bot\setminus \mathfrak G_M^\bot$.
 If, in addition, the system $\varphi(x\dot-h)_{h\in H_0}$ is
 orthonormal, then $\varphi(x)$ generate an orthogonal MRA.
 \end{theorem}
So,  to find a refinable function that generates orthogonal MRA,
we need to take  a function $m_0(\chi)$ that satisfies conditions T1,
T2, T3, (\ref{eq3.5}), construct the function
$$
 \hat\varphi(\chi)=\prod\limits_{k=0}^\infty m_0(\chi{\cal
 A}^{-k})\in \mathfrak D_{-N}(\mathfrak G_M^\bot)
 $$
 and check that the system $\varphi(x\dot-h)_{h\in H_0}$ is
 orthonormal.

 For any zero-dimensional group $G$ the shifts system $(\varphi(x\dot-h))_{h\in H_0}$
 is orthonormal if the condition $|\hat\varphi(\chi)|={\bf 1}_{G_0^\bot}(\chi)$ is valid \cite{14}.
 For Vilenkin group $\mathfrak G$ we can give another condition.

\begin{definition}
Let $N,M\in \mathbb N$. A set $E \subset X$ is called
$(N,M)$-elementary if $E$ is disjoint union of $p^N$ cosets
$$
\mathfrak G_{-N}^\bot\zeta_j=\mathfrak
G_{-N}^\bot\underbrace{r_{-N}^{\alpha_{-N}}r_{-N+1}^{\alpha_{-N+1}}\dots
r_{-1}^{\alpha_{-1}}}_{\xi_j}\underbrace{r_{0}^{\alpha_{0}}\dots
r_{M-1}^{\alpha_{M-1}}}_{\eta_j}=\mathfrak G_{-N}^\bot\xi_j\eta_j,
$$
$j=0,1,...,p^N-1,
j=\alpha_{-N}+\alpha_{-N+1}p+\dots+\alpha_{-1}p^{N-1}$
$(\alpha_{\nu}=\overline{0,p-1})$ such that\\
1) $\bigsqcup\limits_{j=0}^{p^N-1}\mathfrak
G_{-N}^\bot\xi_j=\mathfrak G_{0}^\bot$, $\mathfrak
G_{-N}^\bot\zeta_0=\mathfrak G_{-N}^\bot$,\\
2) for any $l=\overline{0,M+N-1}$ the intersection $(\mathfrak
G_{-N+l+1}^\bot\setminus \mathfrak G_{-N+l}^\bot)\bigcap
E\ne\emptyset$.
\end{definition}
\begin{lemma}
The set $H_0\subset \mathfrak G$ is an orthonormal  system on any
$(N,M)$-elementary set $E\subset X$.
\end{lemma}
{\bf Proof.} Using the definition of $(N,M)$-elementary set we
have
$$
\int\limits_E(\chi,h)\overline{(\chi,g)}\,d\nu(x)=\sum_{j=0}^{p^N-1}\int\limits_{\mathfrak
G_{-N}^\bot\zeta_j}(\chi,h)\overline{(\chi,g)}\,d\nu(x)=
$$
$$
 =\sum\limits_{j=0}^{p^N-1}\int\limits_X{\bf 1}_{\mathfrak
G_{-N}^\bot\zeta_j}(\chi)(\chi,h)\overline{(\chi,g)}\,d\nu(x)=
$$
$$
=\sum\limits_{j=0}^{p^N-1}\int\limits_X{\bf 1}_{\mathfrak
G_{-N}^\bot\zeta_j}(\chi\eta_j)(\chi\eta_j,h)\overline{(\chi\eta_j,g)}\,d\nu(x)=
$$
$$
=\sum\limits_{j=0}^{p^N-1}\int\limits_X{\bf 1}_{\mathfrak
G_{-N}^\bot\xi_j}(\chi)(\chi,h)\overline{(\chi,g)}(\eta_j,h)\overline{(\eta_j,g)}\,d\nu(x).
$$
Since
$$
(\eta_j,h)=(r_0^{\alpha_0}r_1^{\alpha_1}\dots
r_{M-1}^{\alpha_{M-1}},a_{-1}g_{-1}\dot+a_{-2}g_{-2}\dot+\dots\dot+a_{-s}g_{-s})=1,
$$
$$
(\eta_j,g)=(r_0^{\alpha_0}r_1^{\alpha_1}\dots
r_{M-1}^{\alpha_{M-1}},b_{-1}g_{-1}\dot+b_{-2}g_{-2}\dot+\dots\dot+b_{-s}g_{-s})=1,
$$
then
$$
\int\limits_E(\chi,h)\overline{(\chi,g)}\,d\nu(x)=\sum\limits_{j=0}^{p^N-1}\int\limits_{\mathfrak
G_{-N}^\bot\xi_j}(\chi,h)\overline{(\chi,g)}\,d\nu(x)=\int\limits_{\mathfrak
G_{0}^\bot}(\chi,h)\overline{(\chi,g)}\,d\nu(x)=
$$
$=\delta_{h,g}$. {$\square$}
\begin{theorem}
Let $(\mathfrak G,\dot+)$ be a $p$-adic Vilenkin group,
$E\subset\mathfrak G_M^\bot$ -- an $(N,M)$-elementary set. If\
 $|\hat\varphi(\chi)|={\bf 1}_E(\chi)$ on $X$ then the system of
shifts $(\varphi(x\dot-h))_{h\in H_0}$ is an orthonormal system on
$\mathfrak G$.
\end{theorem}
{\bf Proof.} Let $\tilde H_0\subset H_0$ be a finite set. Using
the Plansherel equation we have
$$
\int\limits_{\mathfrak
G}\varphi(x\dot-g)\overline{\varphi(x\dot-g)}\,d\mu(x)=
\int\limits_X|\hat\varphi(\chi)|^2\overline{(\chi,g)}(\chi,h)d\nu(\chi)=
\int\limits_E(\chi,h)\overline{(\chi,g)}d\nu(\chi)=
$$
$$
=\sum_{j=0}^{p^N-1}\int\limits_{\mathfrak
G_{-N}^\bot\zeta_j}(\chi,h)\overline{(\chi,g)}\,d\nu(\chi).
$$
Transform the inner integral
$$
\int\limits_{\mathfrak
G_{-N}^\bot\zeta_j}(\chi,h)\overline{(\chi,g)}\,d\nu(\chi)=\int\limits_X{\bf
1}_{\mathfrak
G_{-N}^\bot\zeta_j}(\chi)(\chi,h)\overline{(\chi,g)}\,d\nu(\chi)=
$$
$$
=\int\limits_X{\bf 1}_{\mathfrak
G_{-N}^\bot\zeta_j}(\chi\eta_j)(\chi\eta_j,h\dot-g)\,d\nu(\chi)=\int\limits_X{\bf
1}_{\mathfrak
G_{-N}^\bot\xi_j}(\chi)(\chi\eta_j,h\dot-g)\,d\nu(\chi)=
$$
$$
=\int\limits_{\mathfrak
G_{-N}^\bot\xi_j}(\chi\eta_j,h\dot-g)\,d\nu(\chi).
$$
Repeating the arguments of lemma 3.4 we obtain
$$
\int\limits_{\mathfrak
G}\varphi(x\dot-h)\overline{\varphi(x\dot-g)}\,d\mu(x)=\delta_{h,g}.\;\;\square
$$
\begin{theorem}[\cite{14}]
 Let  $\varphi(x)\in {\mathfrak
 D}_M(\mathfrak G_{-N})$. A shifts system
 $(\varphi(x\dot-h))_{h\in H_0}$ will be orthonormal if and only
if for any
$\alpha_{-N},\alpha_{-N+1},\dots,\alpha_{-1}=\overline{(0,p-1)}$
 \begin{equation} \label{eq36}
 \sum_{\alpha_{0},\alpha_1,\dots,\alpha_{M-1}=0}^{p-1}|\hat\varphi(\mathfrak G_{-N}^\bot
 r_{-N}^{\alpha_{-N}}\dots r_0^{\alpha_0}\dots
 r_{M-1}^{\alpha_{M-1}})|^2=1.
 \end{equation}
 \end{theorem}
 \begin{lemma}[\cite{14}]
Let $\hat\varphi\in \mathfrak D_{-N}(\mathfrak G_M^\bot)$ be a
solution of the refinement equation
$$
\hat\varphi(\chi)=m_0(\chi)\hat\varphi(\chi{\cal A}^{-1})
$$
and $(\varphi(x\dot-h))_{h\in H_0}$ be an orthonormal system.\\
 Then for any
$\alpha_{-N},\alpha_{-N+1},\dots,\alpha_{-1}=\overline{0,p-1}$
 \begin{equation} \label{eq37}
\sum_{\alpha_0=0}^{p-1}|m_0(\mathfrak G_{-N}^\bot
r_{-N}^{\alpha_{-N}}r_{-N+1}^{\alpha_{-N+1}}\dots
r_{-1}^{\alpha_{-1}}r_{0}^{\alpha_{0}})|^2=1.
 \end{equation}
\end{lemma}
\section{Trees and refinable functions}\label{s4}
In this section we reduce the problem of construction of step
refinable function to construction of some tree.

 We will consider some
special class of refinable functions $\varphi(\chi)$ for which
$|\hat\varphi(\chi)|$ is a characteristic function of a set.
Define this class.
\begin{definition}
A mask $m_0(\chi)$ is called $N$-elementary $(N\in\mathbb N_0)$ if
$m_0(\chi)$ is constant on cosets $\mathfrak G_{-N}^\bot\chi$, its
absolute value $|m_0(\chi)|$ has two values only: 0 and 1, and
$m_0(\mathfrak G_{-N}^\bot)=1$. The refinable function
$\varphi(x)$ with Fourier transform
$$
\hat\varphi(\chi)=\prod\limits_{n=0}^\infty m_0(\chi{\cal A}^{-n})
$$
is called $N$-elementary too. $N$-elementary function $\varphi$
 is called $(N,M)$-elementary if  $\hat\varphi(\chi)\in \mathfrak D_{-N}(\mathfrak
 G_{M}^\bot)$. In this case  we will
 call the Fourier transform $\hat\varphi(\chi)$  $(N,M)$-elementary, also.
\end{definition}

\begin{definition}
Let $\tilde
E=\bigsqcup\limits_{\alpha_{-N},\dots,\alpha_{-1},\alpha_0}\mathfrak
G_{-N}^\bot r_{-N}^{\alpha_{-N}}\dots
r_{-1}^{\alpha_{-1}}r_0^{\alpha_{0}}\subset \mathfrak G_1^\bot$ be
an $(N,1)$-elementary set. We say that the set $\tilde E_X$ is a
periodic extension of $\tilde E$ if
$$
\tilde E_X=\bigcup
\limits_{s=1}^\infty\bigsqcup\limits_{\alpha_1,\dots,\alpha_s=0}^{p-1}\tilde
E r_1^{\alpha_1}r_2^{\alpha_2}\dots r_s^{\alpha_s}.
$$
We say that  $\tilde E$ generates an $(N,M)$ elementary set $E$,
if $\bigcap\limits_{n=0}^\infty \tilde E_X{\cal A}^n=E$.
\end{definition}
Since $\tilde E_X\supset \mathfrak G_{-N}^\bot$ then
$\bigcap\limits_{n=0}^{M+1} \tilde E_X{\cal A}^n=E$ and
$\left(\bigcap\limits_{n=0}^{M+1} \tilde E_X{\cal
A}^n\right)\bigcap(\mathfrak G_{M+1}^\bot\setminus \mathfrak
G_{M}^\bot)=\emptyset$. The converse is also true. Since
$$
\left(\bigcap\limits_{n=0}^{M+1} \tilde E_X{\cal
A}^n\right)\bigcap(\mathfrak G_{M+1}^\bot\setminus \mathfrak
G_{M}^\bot)=\emptyset.
$$
Then we have
$$
\left(\bigcap\limits_{n=0}^{M+2} \tilde E_X{\cal
A}^n\right)\bigcap(\mathfrak G_{M+2}^\bot\setminus \mathfrak
G_{M+1}^\bot)=\tilde E_X\bigcap\left(\bigcap\limits_{n=0}^{M+1}
\tilde E_X{\cal A}^n\bigcap(\mathfrak G_{M+1}^\bot\setminus
\mathfrak G_{M}^\bot)\right){\cal A}=
$$
$$
=\tilde E_X\bigcap\emptyset=\emptyset.
$$

 Let $N$ be a natural number. Denote $V=\{0,1.\dots,p-1\}$ and construct a tree $T(V)$ in the
 following way:\\
 1) The root of this tree and its vertices of level $1,2,\dots,N-1$ are equal to zero.\\
 2) Any path $(\alpha_k\to\alpha_{k+1}\to\dots\to\alpha_{k+N-1})$
 of length $N$ is present in the tree $T(V)$ exactly 1 time.

 Such tree  we will call  $N$-valid.

 For example for $p=3, N=2$ we can construct the tree

\begin{picture}(100,50)

  \put(0,20){\circle{6}}
  \put(-2,18){$0$}
  \put(2,20){\vector(1,0){15}}
  \put(18,18){$0$}
  \put(40,20){\circle{6}}
   \put(22,20){\vector(1,0){15}}
  \put(38,18){$2$}
   \put(42,22){\vector(1,1){10}}
   \put(42,20){\vector(1,0){10}}
   \put(42,18){\vector(1,-1){10}}
   \put(20,20){\circle{6}}
   \put(55,35){\circle{6}}
   \put(53,33){$1$}
   \put(55,20){\circle{6}}
   \put(53,18){$0$}
   \put(55,5){\circle{6}}
   \put(53,3){$2$}
   \put(58,35){\vector(1,0){15}}
   \put(76,35){\circle{6}}
   \put(75,33){$0$}
   \put(58,20){\vector(1,0){15}}
   \put(76,20){\circle{6}}
   \put(75,18){$1$}
   \put(78,22){\vector(1,1){11}}
   \put(92,35){\circle{6}}
   \put(90,33){$1$}
   \put(78,18){\vector(1,-1){11}}
   \put(92,5){\circle{6}}
   \put(90,3){$2$}
    \put(0,0){Figure 1}
   \end{picture}

\noindent
 This tree contains  any edge
$$
(0,0),(0,1),(0,1),(1,0),(1,1),(1,2),(2,0),(2,1),(2,2)
$$
 exactly 1 time and height $T(V)=6$.

Using the tree $T(V)$ we will construct the family of cosets in
the
following way:\\
For any  a path
$$
(\alpha_s\to\alpha_{s-1}\to\dots\to\alpha_{s-N+1}\to\alpha_{s-N}\to\alpha_{s-N-1}\to\dots\to\alpha_{-N+1}\to\alpha_{-N})
$$
in which $\alpha_{s}=\alpha_{s-1}=\dots=\alpha_{s-N+1}=0$.\\
we construct cosets
\begin{equation}                                \label{eq4.1}
G_{-N}^\bot r_{-N}^{\alpha_{-N}}r_{-N+1}^{\alpha_{-N+1}}\dots
r_{0}^{\alpha_{0}}, G_{-N}^\bot
r_{-N}^{\alpha_{-N+1}}r_{-N+1}^{\alpha_{-N+2}}\dots
r_{0}^{\alpha_{1}},\dots,G_{-N}^\bot
r_{-N}^{\alpha_{s-N}}r_{-N+1}^{\alpha_{s-N+1}}\dots
r_{0}^{\alpha_{s}},
\end{equation}
\begin{equation}                                \label{eq4.2}
G_{-N}^\bot r_{-N}^{\alpha_{s-N+1}}\dots r_{-1}^{\alpha_{s}},
G_{-N}^\bot r_{-N}^{\alpha_{s-N+2}}\dots
r_{-2}^{\alpha_{s}},\dots,G_{-N}^\bot r_{-N}^{\alpha_{s}}.
\end{equation}
The union of all such cosets we denote as $\tilde E$. It is clear
that $\tilde E\subset G_1^\bot$.

\begin{definition}
Let $\tilde E_X$ be a periodic extension of $\tilde E$. We say
that the tree $T(V)$ generates a set $E$, if
$E=\bigcap\limits_{n=0}^\infty\tilde E_X{\cal A}^n.$
\end{definition}
\begin{lemma}
Let $T(V)$ be a $N$-valid  tree. Let $E\subset X$ be a set
generated by the tree $T(V)$, $H$ -- height of $T(V)$. Then $E$ is
an $(N,H-2N)$-elementary set.
\end{lemma}
{\bf Proof.} Let us denote
$$
m(\chi)={\bf 1}_{\tilde
E_X}(\chi),\;\;M(\chi)=\prod\limits_{n=0}^\infty m(\chi{\cal
A}^{-n}).
$$
First we note that $M(\chi)={\bf 1}_E(\chi)$. Indeed
$$
{\bf 1}_E(\chi)=1\Leftrightarrow\chi\in E\Leftrightarrow
\forall\,n,\;\chi{\cal A}^{-n}\in\tilde
E_X\Leftrightarrow\forall\,n,\;{\bf 1}_{\tilde E_X}(\chi{\cal
A}^{-n})=1\Leftrightarrow
$$
$$
\forall\,n,\;m(\chi{\cal
A}^{-n})=1\Leftrightarrow\prod\limits_{n=0}^\infty m(\chi{\cal
A}^{-n})=1\Leftrightarrow M(\chi)=1.
$$
It means that $M(\chi)={\bf 1}_E(\chi)$.\\
Now we will prove, that ${\bf 1}_E(\chi)=0$ for $\chi\in\mathfrak
G_{H-2N+1}^\bot\setminus\mathfrak G_{H-2N}^\bot$. Since $\tilde
E_X\supset \mathfrak G_{-N}^\bot$ it follows that ${\bf 1}_{\tilde
E_X}(\mathfrak G_{H-2N}^\bot{\cal A}^{-H+N})={\bf 1}_{\tilde
E_X}(\mathfrak G_{-N}^\bot)=1$. Consequently
$$
\prod\limits_{n=0}^\infty{\bf 1}_{\tilde E_X}(\chi{\cal
A}^{-n})=\prod\limits_{n=0}^{H-N-1}{\bf 1}_{\tilde E_X}(\chi{\cal
A}^{-n})
$$
if $\chi\in \mathfrak G_{H-2N+1}^\bot\setminus\mathfrak
G_{H-2N}^\bot$. Let us denote
$$
m(\mathfrak G_{-N}^\bot
r_{-N}^{\alpha_{-N}}r_{-N+1}^{\alpha_{-N+1}}\dots
r_0^{\alpha_0})=\lambda_{\alpha_{-N},\alpha_{-N+1},\dots,\alpha_0}.
$$
By the definition of cosets \eqref{eq4.1}, \eqref{eq4.2} $
m(\mathfrak G_{-N}^\bot
r_{-N}^{\alpha_{-N}}r_{-N+1}^{\alpha_{-N+1}}\dots r_0^{\alpha_0})
\ne 0\Leftrightarrow$ the vector
$(\alpha_0,\alpha_1,...,\alpha_{-N+1},\alpha_{-N})$ is a path
$(\alpha_0\rightarrow\alpha_1\rightarrow...\rightarrow\alpha_{-N+1}\rightarrow\alpha_{-N})$
 of
the tree $T(V)$.

We need to prove that
$$
{\bf 1}_E(\mathfrak G_{-N}^\bot
r_{-N}^{\alpha_{-N}}r_{-N+1}^{\alpha_{-N+1}}\dots
r_{H-2N}^{\alpha_{H-2N}})=0
$$
for $\alpha_{H-2N}\ne 0$. Since $\tilde E_X$ is a periodic
extension of $\tilde E$ it follows that the function $m(\chi)={\bf
1}_{\tilde E_X}(\chi)$ is periodic with any period
$r_{1}^{\alpha_{1}}r_{2}^{\alpha_{2}}\dots r_{s}^{\alpha_{s}}$,
$s\in\mathbb N$, i.e. $m(\chi
r_{1}^{\alpha_{1}}r_{2}^{\alpha_{2}}\dots
r_{s}^{\alpha_{s}})=m(\chi)$ when $\chi\in \mathfrak G_1^\bot$.
Using this fact we can write $M(\chi)$ for  $\chi\in \mathfrak
G_{H-2N+1}^\bot\setminus \mathfrak G_{H-2N}^\bot$ in the form
$$
M(\mathfrak G_{-N}^\bot\zeta)=M(\mathfrak G_{-N}^\bot
r_{-N}^{\alpha_{-N}}r_{-N+1}^{\alpha_{-N+1}}\dots
r_{H-2N}^{\alpha_{H-2N}})=
$$
$$
 =m(\mathfrak G_{-N}^\bot r_{-N}^{\alpha_{-N}}
 r_{-N+1}^{\alpha_{-N+1}}\dots r_{0}^{\alpha_{0}})
 m(\mathfrak  G_{-N}^\bot r_{-N}^{\alpha_{-N+1}}r_{-N+1}^{\alpha_{-N+2}}\dots
 r_{0}^{\alpha_{1}})\dots
$$
\begin{equation}  \label{eq4.3}
 m(\mathfrak  G_{-N}^\bot
r_{-N}^{\alpha_{H-3N}}r_{-N+1}^{\alpha_{H-3N+1}}\dots
r_{-1}^{\alpha_{H-2N-1}}
 r_{0}^{\alpha_{H-2N}})
 \end{equation}
$$
m(\mathfrak  G_{-N}^\bot
r_{-N}^{\alpha_{H-3N+1}}r_{-N+1}^{\alpha_{H-3N+1}}\dots
r_{-1}^{\alpha_{H-2N}})\dots m(\mathfrak G_{-N}^\bot r_{-N}^{\alpha_{H-2N-1}}r_{-N+1}^{\alpha_{H-2N}})m(\mathfrak G_{-N}^\bot r_{-N}^{\alpha_{H-2N}}).
$$
Assume that $M(\mathfrak G_{-N}^\bot\zeta)\neq 0$. Then all
factors in \eqref{eq4.3} are nonzero. So we have the path
$$
0\rightarrow\dots\rightarrow 0\rightarrow\alpha _{H-2N}\neq 0\rightarrow \alpha _{H-2N-1}\rightarrow\dots\rightarrow \alpha
_0\rightarrow \dots \rightarrow\alpha _{-N+1} \rightarrow\alpha
_{-N},
$$
where there are $N$ zeroes at the beginning of the path. The length of such path is $H+1$, which contradicts the condition $height(T)=H$.

Now we prove that $E$ is $(1,H-2N)$ elementary set.
Since the tree T(V) is $N$-valid, it has all possible combinations of $N$ numbers $\alpha_i=\overline{0,p-1}$ as its paths, and we have the first property of elementary sets satisfied. Also, since $height(T)=H$, there exists a path
$$\alpha_1=0\rightarrow\dots\rightarrow \alpha_N=0\rightarrow\alpha_{N+1}\neq 0\rightarrow\alpha_{N+2}\rightarrow\dots\rightarrow\alpha_{H}$$
of length $H$. Such path generates the set $\mathfrak G_{-N}^\bot r_{-N}^{\alpha_{N+1}}\subset\mathfrak G_{-N+1}\setminus\mathfrak G_{-N}$, since $\alpha_{N+1}\neq 0$. Also, the same path generates the set $\mathfrak G_{-N}^\bot r_{-N}^{\alpha_{N+2}}r_{-N+1}^{\alpha_{N+1}}\subset\mathfrak G_{-N+2}\setminus\mathfrak G_{-N+1}$. Continuing this process we will obtain all sets $$\forall l=\overline{0,H-N-1}, \mathfrak G_{-N}^\bot\prod\limits_{n=0}^l r_{-N+n}^{\alpha_{N+1+n}}\subset\mathfrak G_{-N+l+1}\setminus\mathfrak G_{-N+l},$$ which means the second property of elementary sets is also satisfied. Thus we can conclude that $E$ is $(1,H-2N)$-elementary set and the lemma is proved.
  $\square$.
\large
 \begin{theorem}
Let $M,p\in\mathbb N$, $p\ge 3$. Let $E\subset \mathfrak G_M^\bot$
be an $(N,M)$-elementary set, $\hat\varphi\in\mathfrak
D_{-N}(\mathfrak G_M^\bot)$, $|\hat\varphi(\chi)|={\bf
1}_E(\chi)$, $\hat\varphi(\chi)$ the solution of the equation
\begin{equation}              \label{eq4.4}
\hat\varphi(\chi)=m_0(\chi)\hat\varphi(\chi{\cal A}^{-1}),
\end{equation}
where $m_0(\chi)$ is an N-elementary mask. Then there exists a
rooted tree $T(V)$ with ${\rm height}(T)=M+2N$ that generates the
set $E$.
\end{theorem}
 {\bf Prof.} Since the set $E$ is $(N,M)$-elementary set and
 $|\hat\varphi(\chi)|={\bf 1}_E(\chi)$, it follows from theorem 3.2
 that the system $(\varphi(x\dot-h))_{h\in H_0}$ is an orthonormal
 system in $L_2(\mathfrak G)$. Using the theorem 3.3 we obtain that
 $\forall\alpha_{-N},\dots,\alpha_{-1}=\overline{0,p-1}$
 $$
 \sum_{\alpha_0,\alpha_1,\dots,\alpha_{M-1}=0}^{p-1}|\hat\varphi(\mathfrak
 G_{-N}^\bot r_{-N}^{\alpha_{-N}}\dots
 r_{-1}^{\alpha_{-1}}r_{0}^{\alpha_0}\dots
 r_{M-1}^{\alpha_{M-1}})|^2=1.
 $$
 Since $\hat\varphi$ is a solution of refinement equation
 \eqref{eq4.4} it follows from lemma 3.5 that
 $\forall\alpha_{-N},\dots,\alpha_{-1}=\overline{0,p-1}$
 \begin{equation}              \label{eq4.5}
 \sum_{\alpha_0=0}^{p-1}|m_0(\mathfrak G_{-N}^\bot r_{-N}^{\alpha_{-N}}
 r_{-1}^{\alpha_{-1}}r_{0}^{\alpha_{0}})|^2=1.
 \end{equation}
 Let as denote $\lambda_{\alpha_{-N},\dots,\alpha_{-1},\alpha_0}:=m_0(\mathfrak
 G_{-N}^\bot r_{-N}^{\alpha_{-N}}\dots r_{-1}^{\alpha_{-1}}r_{0}^{\alpha_0})$. Then we
 write \eqref{eq4.5} in the form
 \begin{equation}              \label{eq4.6}
 \sum_{\alpha_0=0}^{p-1}|\lambda_{\alpha_{-N},\dots,\alpha_{-1},\alpha_0}|^2=1.
 \end{equation}
 Since the mask $m_0(\chi)$ is N-elementary it follows that
 $|\lambda_{\alpha_{-N},\dots,\alpha_{-1},\alpha_0}|$ takes one of two values only: 0 or 1.

Now we will construct the tree $T$. We will begin with the path of $N$ zeros
$$0_1\rightarrow 0_2\rightarrow\dots\rightarrow 0_N,$$
where $0_1$ is the root of the tree.

 Let $\mathfrak U$ be a family
of cosets $\mathfrak G_{-N}^\bot\zeta\subset \mathfrak G_M^\bot$
such that $\hat\varphi(\mathfrak G_{-N}^\bot \zeta)\ne 0$ and
$\mathfrak G_{-N}^\bot\notin \mathfrak U$. We can write a coset
$\mathfrak G_{-N}^\bot\zeta\in \mathfrak U$ in the form
 \begin{equation}              \label{eq4.7}
 \mathfrak G_{-N}^\bot\zeta=
\mathfrak G_{-N}^\bot r_{-N}^{\alpha_{-N}}\dots
r_{-1}^{\alpha_{-1}}r_{0}^{\alpha_{0}}\dots r_{s}^{\alpha_{s}},\
\alpha_{s}\neq 0.
\end{equation}
 Here $s\leq M-1$ since each coset in
 $\mathfrak U$ is a subset of $\mathfrak G_{M}^\bot$,
 and there exists at least one coset with $s=M-1$ since function is
 (N,M)-elementary. If $s=M-1$ and $\alpha_{s+1}+\dots +\alpha_{s+l}\neq 0$ then  coset
 $$
   \mathfrak G_{-N}^\bot r_{-N}^{\alpha_{-N}}\dots
r_{-1}^{\alpha_{-1}}r_{0}^{\alpha_{0}}\dots
r_{s}^{\alpha_{s}}r_{s+1}^{\alpha_{s+1}}\dots
r_{s+l}^{\alpha_{s+l}}\notin \mathfrak U
 $$

  0)Initially, we take a coset
 $$
    \mathfrak G_{-N}^\bot\zeta_1=
\mathfrak G_{-N}^\bot r_{-N}^{\alpha_{-N}^{(1)}}\dots
r_{-1}^{\alpha_{-1}^{(1)}}r_{0}^{\alpha_{0}^{(1)}}\dots
r_{s_1}^{\alpha_{s_1}^{(1)}}\in  \mathfrak U,\  \alpha_{s_1}\neq 0
 $$
 and  connect the path
 $$
 p^{(1)}= \alpha_{s_1}^{(1)}\rightarrow\dots\rightarrow\alpha_0^{(1)}\rightarrow\alpha_{-1}^{(1)}\rightarrow\dots\rightarrow\alpha_{-N}^{(1)}
 $$
 to the $0_N$ vertex. We obtain the tree $T^{(0)}$ that contains  unique branch
$$
 T^{(0)}=(0_1\rightarrow 0_2\rightarrow\dots\rightarrow 0_N \rightarrow
  \alpha_{s_1}\rightarrow\dots\rightarrow\alpha_0\rightarrow\alpha_{-1}\rightarrow\dots
  \rightarrow\alpha_{-N}).
 $$
  1) On the {\bf first step}, take  another coset
 $$
 \mathfrak G^\bot_{-N}\zeta_2 = \mathfrak G_{-N}^\bot r_{-N}^{\alpha_{-N}^{(2)}}\dots r_{-1}^{\alpha_{-1}^{(2)}}r_{0}^{\alpha_{0}^{(2)}}\dots
  r_{s_2}^{\alpha_{s_2}^{(2)}}\in \mathfrak U\setminus \mathfrak G_{-N}^\bot\zeta_1,\alpha_{s_2}^{(2)}\neq 0
 $$
  which generates the path
 $$
  p^{(2)}= (
  \alpha_{s_2}^{(2)}\rightarrow\dots\rightarrow\alpha_0^{(2)}
  \rightarrow\alpha_{-1}^{(2)}\rightarrow\dots\rightarrow\alpha_{-N}^{(2)})
 $$

  Let us add the path  $0_1\rightarrow 0_2\rightarrow\dots\rightarrow 0_N$ to
 the beginning of the path  $p^{(2)}$ and denote it as $\tilde{p}^{(2)}$, i.e.
 $$
  \tilde{p}^{(2)}=(0_1\rightarrow 0_2\rightarrow\dots\rightarrow 0_N \rightarrow
  \alpha_{s_2}^{(2)}\rightarrow\dots\rightarrow\alpha_0^{(2)}\rightarrow\alpha_{-1}^{(2)}
  \rightarrow\dots\rightarrow\alpha_{-N}^{(2)})
 $$
 Now we will
 include   this path into our tree $T^{(0)}$. To include it we will compare the path
 $\tilde{p}^{(2)}$ with the tree  $T^{(0)}$. \\

 There are 3 possible cases:\\
 1)The path $p^{(0)}$ is shorter than $p^{(1)}$ and\\
 $$
 \alpha_{s_0}^{(0)}=\alpha_{s_1}^{(1)},\alpha_{s_0-1}^{(0)}=\alpha_{s_1-1}^{(1)},...,
 \alpha_{-N}^{(0)}=\alpha_{s_1-(s_0+N)}^{(1)}.
 $$
 In this case we connect  the tail
 $$
  \alpha_{s_1-(s_0+N+1)}^{(1)} \rightarrow
  \alpha_{s_1-(s_0+N+2)}^{(1)} \rightarrow \dots \rightarrow \alpha_{-N}^{(1)}
 $$
 of the path $p^{(1)}$ to the vertex $\alpha_{-N}^{(0)}$.\\
 2)The path $p^{(0)}$ is longer than $p^{(1)}$ and\\
 $$
 \alpha_{s_0}^{(0)}=\alpha_{s_1}^{(1)},\alpha_{s_0-1}^{(0)}=\alpha_{s_1-1}^{(1)},...,
  \alpha_{s_0-(s_1+N)}^{(0)}=\alpha_{-N}^{(1)}.
 $$
 In this case the path $\tilde{p}^{(1)}$ is already a path of the tree
 $T^{(0)}$ and we leave the tree
 $T^{(0)}$ unchanged.\\
 3)There exists an integer $l$ such that
 $\alpha_{s_1-l}^{(1)}\neq \alpha_{s_0-l}^{(0)}$ and
 $\forall k<l,\alpha_{s_1-k}^{(1)}= \alpha_{s_0-k}^{(0)}$. If $l=-1$ then
 we get $\alpha_{s_1-l}^{(1)}=0_N$. When $l$ is calculated we connect the path \\
 $$
 \alpha_{s_1-l}^{(1)} \rightarrow \alpha_{s_1-l-1}^{(1)}
 \rightarrow\dots \rightarrow \alpha_{-N}^{(1)}
 $$
 to the vertex $\alpha_{s_1-l+1}^{(0)}$ and obtain the tree

 
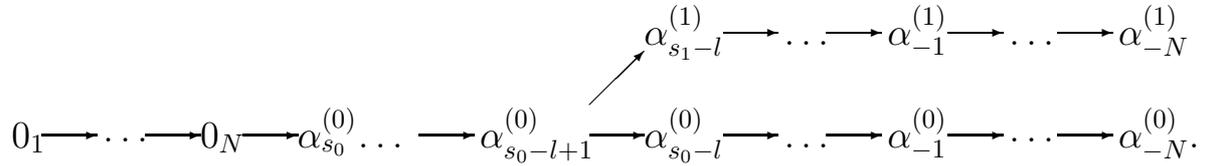
\begin{picture}(160,35)
  \put(-7,9){$0_1$}
  \put(-2,11){\vector(1,0){9}}

    \put(8,10){$\dots$}
   \put(15,11){\vector(1,0){9}}

   \put(24,9){$0_N$}
   \put(31,11){\vector(1,0){9}}

   \put(40,9){$\alpha_{s_0}^{(0)}$}

    \put(50,9){$\dots $}
    \put(60,11){\vector(1,0){9}}
    \put(70,9){$\alpha_{s_0-l+1}^{(0)}$}
   \put(88,11){\vector(1,0){9}}

    \put(97,9){$\alpha_{s_0-l}^{(0)}$}
    \put(110,11){\vector(1,0){9}}

    \put(120,9){$\dots$}


    \put(127,11){\vector(1,0){9}}
    \put(137,9){$\alpha_{-1}^{(0)}$}

     \put(147,11){\vector(1,0){9}}

    \put(157,10){$\dots$}
    \put(165,11){\vector(1,0){9}}
    \put(175,9){$\alpha_{-N}^{(0)}.$}

    \put(88,16){\vector(1,1){9}}
    \put(97,26){$\alpha_{s_1-l}^{(1)}$}

    \put(110,28){\vector(1,0){9}}

    \put(120,26){$\dots$}

    \put(127,28){\vector(1,0){9}}
    \put(137,26){$\alpha_{-1}^{(1)}$}

     \put(147,28){\vector(1,0){9}}

    \put(157,26){$\dots$}
    \put(165,28){\vector(1,0){9}}
    \put(175,26){$\alpha_{-N}^{(1)}$}

    \put(70,0){Figure 2}
   \end{picture}

 This is the end of first step.

   Consider $n$ steps fulfilled, i.e. paths  $p^{(0)},p^{(1)},...,p^{(n)}$
   are chosen and the correspondent tree $T^{(n)}$ is constructed. Now we
   will perform the $(n+1)$-th step.
   Let us take a coset
  $$
  \mathfrak G_{-N}^\bot\zeta_{n+1}=\mathfrak G_{-N}^\bot r_{-N}^{\alpha_{-N}^{(n+1)}}\dots
  r_{-1}^{\alpha_{-1}^{(n+1)}}r_{0}^{\alpha_{0}^{(n+1)}}\dots
  r_{s_{n+1}}^{\alpha_{s_{n+1}}^{(n+1)}}\in  \mathfrak U\setminus\bigcup\limits_{k=1}^n \mathfrak G_{-N}^\bot\zeta_k,\  \alpha_{s_{n+1}}^{(n+1)}\neq 0,
  $$
  which generates a path
  $$
  p^{(n+1)}=(\alpha_{s_{n+1}}^{(n+1)}\rightarrow\dots\rightarrow
  \alpha_0^{(n+1)}\rightarrow\alpha_{-1}^{(n+1)}\rightarrow\dots\rightarrow\alpha_{-N}^{(n+1)}).
  $$
  and denote
  $$
  \tilde{p}^{(n+1)}=(0_1\rightarrow \dots\rightarrow 0_N \rightarrow
  \alpha_{s_{n+1}}^{(n+1)}\rightarrow\dots\rightarrow\alpha_0^{(n+1)}
  \rightarrow\alpha_{-1}^{(n+1)}\rightarrow\dots\rightarrow\alpha_{-N}^{(n+1)}).
  $$
  Now we will include the path $\tilde{p}^{(n+1)}$ into the tree
  $T^{(n)}$. To do it, we will be looking for a path in the tree $T^{(n)}$ such that it has the longest starting sequence matching with the beginning of $\tilde{p}^{(n+1)}$. \\
   {\bf Step $n+1.1.$} If $\alpha_{s_{n+1}}^{(n+1)}$ is not equal to any
   vertex of level $N+1$ of the tree $T^{(n)}$ then we connect the
   path $p^{(n+1)}$ to the vertex $0_N$, obtain the new tree
   $T^{(n)}$ and finish the step.\\
   {\bf Step $n+1.2.$}
   Otherwise there exists such (always unique) vertex of the level $N+1$, which we will denote by $\alpha_{(N+1),i}$, equal to $\alpha_{s_{n+1}}^{(n+1)}$ we consider all vertices of level $N+2$ connected to it. If there are no vertices connected or there are no such vertices matching $\alpha_{s_{n+1}-1}^{(n+1)}$ then we connect the tail of $p^{(n+1)}$ starting from the element $\alpha_{s_{n+1}-1}^{(n+1)}$ to the vertex $\alpha_{(N+1),i_j}$, obtain new tree and finish the step. Otherwise, if there exists vertex of level $N+2$ $\alpha_{(N+2),i}$ equal to $\alpha_{s_{n+1}-1}^{(n+1)}$, we continue the process of including the path $p^{(n+1)}$ into the tree until either there are no more elements in the path $p^{(n+1)}$ or at some level there are no vertices equal to corresponding element of the path $p^{(n+1)}$. In the first case the tree is left unchanged at this step. In the second case the tail of $p^{(n+1)}$ is added to the tree somewhere. Obviously, since the path $p^{(n+1)}$ has finite number of elements the process will also be finite.

The description of the $(n+1)$-th step is finished and there are only few final remarks left.

1)The resulting graph is a tree, since we produce no cycles at each step.

2)The process of constructing such tree is finite, i.e. contains finite number of steps since during each step we use different coset of $\mathfrak U$ and there is a finite number of such cosets.

So, at this point we have obtained a tree. Let us prove that this tree $T$ is N-valid. To prove it, we must show, that each path of $N$ elements is unique in our tree. Firstly, let us prove that the path of $N$ zeros appears only once in our tree -- and it is the path starting from its root. Indeed, let us assume that the path exists somewhere else in the tree $T$ and that it is a part of some path
$$0_1\rightarrow\dots\rightarrow 0_N\rightarrow\alpha_s\rightarrow\dots\rightarrow\alpha_k\rightarrow 0_1\rightarrow\dots\rightarrow 0_N\rightarrow\dots\rightarrow\alpha_{-N}$$
from root to leaf. Since $\alpha_s\neq 0$ there exists at least one nonzero element between two instances of the path $0_1\rightarrow\dots\rightarrow 0_N$. Without the loss of generality we can consider $\alpha_k\neq 0$.

Using the same technique as in \eqref{eq4.3}, we can conclude, that
$$|\hat\varphi(\mathfrak G_{-N}^\bot r_{-N}^{\alpha_{-N}}\dots r_{k-2N-1}^{0_{N}}\dots r_{k-N-1}^{0_{1}}r_{k-N}^{\alpha_k}\dots r_s^{\alpha_s})|=$$
$$=|\lambda_{\alpha_{-N},\dots,\alpha_{-1},\alpha_0}|\dots|\lambda_{0_{N},\dots,0_{1},\alpha_k}|\dots|\lambda_{\alpha_s,0,\dots,0}|=1,$$
which in particular means that $|\lambda_{0,\dots,0,\alpha_k}|=1$. Also, by the properties of the mask $\lambda_{0,\dots,0,0}=1$. These equalities contradict \eqref{eq4.6}.

Now, let us assume that the arbitrary path $\gamma_{-1}\rightarrow\dots\rightarrow\gamma_{-N}$ appears twice. Thus, it is a subpath of 2 different paths from root to leaf
 $$0_1\rightarrow\dots 0_N\rightarrow\alpha_s\rightarrow\dots\rightarrow\alpha_{k}\rightarrow\gamma_{-1}\rightarrow\dots\rightarrow\gamma_{-N}\rightarrow\dots\rightarrow\alpha_{-N},\ k<s,$$
$$0_1\rightarrow\dots 0_N\rightarrow\beta_{s'}\rightarrow\dots\rightarrow\beta_{k'}\rightarrow\gamma_{-1}\rightarrow\dots\rightarrow\gamma_{-N}\rightarrow\dots\rightarrow\beta_{-N},\  k'<s'.$$
Let us denote $0_i=\alpha_{s+N-i+1}=\beta_{s'+N-i+1}$. Now, let us prove, that $\exists j\geqslant 0:\alpha_{k+j}\neq\beta_{k'+j}$.

We assume that the length of $\alpha$ subpath is less than the length of $\beta$ subpath, i.e. $s-k<s'-k'$. Firstly, let us check if $\alpha_k\neq\beta_{k'}$. If they are equal, let's check if $\alpha_{k+1}\neq\beta_{k'+1}$. If we haven't encountered nonequal pair before $0_N=\alpha_{s+1}$ and $\beta_{k'-k+s+1}$, we check if they are nonequal. If not (i.e they are equal), we check all the remaining pairs. If next $N-1$ elements of $\beta$ subpath are equal to elements $0_i$ of the $\alpha$ subpath, it contradicts the fact that there is only one subpath of $N$ zeros in our tree. Thus in this case $\exists j\geqslant 0:\alpha_{k+j}\neq\beta_{k'+j}$.

Now, let us assume, that both subpaths are of the same length. If $\forall j\geqslant 0:\alpha_{k+j}=\beta_{k'+j}$ then, by construction of the tree $T$ these two paths actually correspond to the same vertices from $0_1$ to $\gamma_{-N}$, which means subpath $\gamma$ does not appear twice in our tree. It contradicts the initial assumption that it does appear twice. Thus in this case $\exists j\geqslant 0:\alpha_{k+j}\neq\beta_{k'+j}$, too.

Let us assume, without loss of generality, that $\alpha_k\neq\beta_{k'}$. Using the same technique as in \eqref{eq4.3}, we can conclude, that
$$|\hat\varphi(\mathfrak G_{-N}^\bot r_{-N}^{\alpha_{-N}}\dots r_{k-2N-1}^{\gamma_{-N}}\dots r_{k-N-1}^{\gamma_{-1}}r_{k-N}^{\alpha_k}\dots r_s^{\alpha_s})|=$$
$$=|\lambda_{\alpha_{-N},\dots,\alpha_{-1},\alpha_0}|\dots|\lambda_{\gamma_{-N},\dots,\gamma_{-1},\alpha_k}|\dots|\lambda_{\alpha_s,0,\dots,0}|=1,$$
$$|\hat\varphi(\mathfrak G_{-N}^\bot r_{-N}^{\beta_{-N}}\dots r_{k'-2N-1}^{\gamma_{-N}}\dots r_{k'-N-1}^{\gamma_{-1}}r_{k'-N}^{\beta_{k'}}\dots r_{s'}^{\beta_{s'}})|=$$
$$=|\lambda_{\beta_{-N},\dots,\beta_{-1},\beta_0}|\dots|\lambda_{\gamma_{-N},\dots,\gamma_{-1},\beta_{k'}}|\dots|\lambda_{\beta_{s'},0,\dots,0}|=1.$$

That means, in particular, that $|\lambda_{\gamma_{-N},\dots,\gamma_{-1},\beta_{k'}}|=|\lambda_{\gamma_{-N},\dots,\gamma_{-1},\alpha_{k}}|=1$, which contradicts \eqref{eq4.6}. Thus our tree is N-valid.

It is evident that this tree generates refinable function
$\hat\varphi$ with a mask $m_0$. Let's show that ${\rm height}(T)=M+2N$.
Indeed, since $\hat\varphi\in  \mathfrak D_{-N}(\mathfrak
 G_{M}^\bot)$ it follows that there exists a coset $\mathfrak G_{-N}^\bot r_{-N}^{\alpha_{-N}}\dots r_{-1}^{\alpha_{-1}}r_{0}^{\alpha_{0}}
 \dots r_{M-1}^{\alpha_{M-1}}$, $\alpha_{M-1}\ne 0$ for which $|\hat\varphi(\mathfrak G_{-N}^\bot r_{-N}^{\alpha_{-N}}\dots r_{-1}^{\alpha_{-1}}r_{0}^{\alpha_{0}}
 \dots r_{M-1}^{\alpha_{M-1}})|=1$. This coset generates a path
 $$0_1\rightarrow\dots\rightarrow 0_N\rightarrow\alpha_{M-1}\rightarrow\dots\rightarrow\alpha_0\rightarrow\alpha_{-1}\rightarrow\dots\rightarrow\alpha_{-N}$$
 of $T$. This path contain $M+2N$ vertices. It
 means that ${\rm height}(T)\ge M+2N$. On the other hand there is no
 coset $\mathfrak G_{-N}^\bot\zeta\subset\mathfrak G_{M+1}^\bot\setminus \mathfrak
 G_{M}^\bot$, consequently  there is no path with $L>M+2N$. So
${\rm height}(T)=M+2N$. The theorem is
proved.  $\square$

\begin{definition}
 Let $T(V)$ be an N-valid tree,  $H=height(T)$. Using cosets \eqref{eq4.1} we define
 the mask $m_0(\chi)$ in the subgroup $\mathfrak G_{1}^\bot$ as
 follows: $m_0(\mathfrak G_{-N}^\bot)=1, m_0(\mathfrak G_{-N}^\bot r_{-N}^{\alpha_{-N}}\dots
 r_{-1}^{\alpha_{-1}}r_0^{\alpha_0})=
 \lambda_{\alpha_{-N},\dots,\alpha_{-1},\alpha_{0}}$,
 $|\lambda_{\alpha_{-N},\dots,\alpha_{-1},\alpha_{0}}|=1$ when
 $\mathfrak G_{-N}^\bot r_{-N}^{\alpha_{-N}}\dots
 r_{-1}^{\alpha_{-1}}r_0^{\alpha_0}\subset \tilde E$,
 $m_0(\mathfrak G_{-N}^\bot r_{-N}^{\alpha_{-N}}\dots
 r_{-1}^{\alpha_{-1}}r_0^{\alpha_0})=
 \lambda_{\alpha_{-N},\dots,\alpha_{-1},\alpha_{0}}=0$ when
 $\mathfrak G_{-N}^\bot r_{-N}^{\alpha_{-N}}\dots
 r_{-1}^{\alpha_{-1}}r_0^{\alpha_0}\subset \mathfrak G_1^\bot\backslash\tilde
 E$.
Let us extend the mask $m_0(\chi)$ on the $X\setminus \mathfrak
G_{1}^\bot$ periodically, i.e. $m_0(\chi
r_1^{\alpha_1}r_2^{\alpha_2}\dots r_s^{\alpha_s})=m_0(\chi)$. Then
we say that the tree $T(V)$ generates the mask $m_0(\chi)$. Set
$\hat \varphi(\chi)=\prod\limits_{n=0}^\infty m_0(\chi{\cal
A}^{-n})$.
It follows from lemma 4.1 that\\
 1) ${\rm
supp}\,\hat\varphi(\chi)\subset\mathfrak G_{H-2N}^\bot$,\\
 2)
$\hat\varphi(\chi)$ is $(N,H-2N)$ elementary function,\\
 3)
$(\varphi(x\dot-h))_{h\in H_0}$ is an orthonormal system.\\ In
this case we say that the tree $T(V)$ generates the refinable
function $\varphi(x)$.
\end{definition}

\begin{theorem}
Let $p\ge 3$ be a prime number,
$T(V)$ an N-valid tree. Let $H$ be are height
of $T(V)$. By $\varphi(x)$ denote the function generated by the
$T(V)$. Then $\varphi(x)$ generates an orthogonal MRA on $p$-adic
Vilenkin group.
\end{theorem}
{\bf Proof.} Since  $T(V)$ generates the the function  $\varphi$
 then 1)$\hat\varphi\in\mathfrak
D_{-N}(\mathfrak G_M^\bot)$, 2)$\hat\varphi(\chi)$ is $(N,H-2N)$-
elementary function, 3)$\hat\varphi(\chi)$ is a solution of
refinable equation \eqref{eq3.3}, 4)$(\varphi(x\dot-h))_{h\in
H_0}$ is an orthonormal system. From the theorem 3.1 it follows
that $\varphi(x)$ generates an orthogonal MRA. $\square$


\section{Construction of wavelet bases}\label{s5}
In \cite{6} and \cite{7} Yu.A.Farkov reduces the problem of
$p$-wavelet decomposition into a problem of matrix extension. We
will use another  method \cite{13}.

As usual, $W_n$ stands for the orthogonal complement of $V_n$ in
$V_{n+1}$: that is $V_{n+1}=V_n\oplus W_n$ and $V_n\bot W_n$
($n\in\mathbb Z$, and $\oplus$ denotes the direct sum).\\
It is readily seen that\\
1) $f\in W_n\Leftrightarrow f({\cal A}x)\in W_{n+1}$,\\
2) $W_n\bot W_k$ for $k\ne n$,\\
3) $\oplus W_n=L_2(\mathfrak G)$, $n\in\mathbb Z$.

From theorems 4.1, 4.2 we derive an algorithm for constructing
wavelet bases.\\
{\bf Step 1.} Choose an arbitrary tree
$T$ - $N$-valid. Let $H$ be a height of the tree $T$.\\
{\bf Step 2.} Choose a finite sequence
$(\lambda_{\alpha_{-N},\dots,\alpha_0})_{\alpha_{-N},\dots,\alpha_0=0}^{p-1}$ such that $\lambda_{0,0,\dots,0}=1,$
$|\lambda_{\alpha_{-N},\dots,\alpha_0}|=1$ if there exists subpath
$\alpha_{-N}\rightarrow\dots\rightarrow\alpha_0$
in the tree $T$, $|\lambda_{\alpha_{-N},\dots,\alpha_0}|=0$ otherwise.\\
{\bf Step 3.} Construct the mask $m_0(\chi)$ and Fourier transform
$\hat\varphi(\chi)$ using definition 4.4. It is clear
that $E={\rm supp}(\hat\varphi(\chi))$ is $(N,H-2N)$-elementary set.\\
{\bf Step 4.} Find coefficients $\beta_h$ for which
\begin{equation}              \label{eq5.1}
m_0(\chi)=\frac{1}{p}\sum_{h\in
H_0^{(N+1)}}\beta_h\overline{(\chi{\cal A}^{-1},h)}.
\end{equation}
To find coefficients $\beta_h$, we write this equation in the form
\begin{equation}\label{eq5.2}
m_0(\chi_k)=\frac{1}{p}\sum_{j=0}^{p^{N+1}-1}\beta_j\overline{(\chi_k,{\cal
A}^{-1}h_j)}
\end{equation}
 where
$$
\begin{array}{ll}
  h_j=a_{-1}g_{-1}\dot+a_{-2}g_{-2}\dot+\dots\dot+a_{-N-1}g_{-N-1}, & \chi_k\in \mathfrak G_{-N}^\bot
r_{-N}^{\alpha_{-N}}\dots r_{-1}^{\alpha_{-1}}r_{0}^{\alpha_{0}}, \\
    j=a_{-1}+a_{-2}p+\dots +a_{-N-1}p^{N},& k=\alpha_{-N}+\dots+\alpha_{-1}p^{N-1}+\alpha_{0}p^{N},
  \\
 a_{-1},a_{-2},\dots,a_{-N}=\overline{0,p-1}, & \alpha_{-N-1},\dots,\alpha_{-1},\alpha_0=\overline{0,p-1}.\\
\end{array}
$$

Since the matrix $\frac{1}{p}\overline{(\chi_k,{\cal A}^{-1}h_j)}$
of this system is unitary it follows that the system \eqref{eq5.2}
has a unique solution.\\
 {\bf Step 5.} We set $m_l(\chi)=m_0(\chi r_0^{-l})$,
 $l=\overline{1,p-1}$, $X_0=\{\chi:|m_0(\chi)|=1\}$. Clearly, $m_l(\chi)$ may be written as
 $$
  m_l(\chi)=\frac{1}{p}\sum_{h\in H_0^{(N+1)}}\beta_h\overline{(\chi
  r_0^{-l},{\cal A}^{-1}h)}=\frac{1}{p}\sum_{h\in
  H_0^{(N+1)}}\beta_h^{(l)}\overline{(\chi,{\cal A}^{-1}h)}
 $$
 where $\beta_h^{(l)}=\beta_h(r_0^l,{\cal A}^{-1}h)$. By the
 construction of $m_l(\chi)$ we have $|m_l(X_0 r_0^l)|=1$.
 From the necessary condition \eqref{eq37} it follows that
 $|m_l(X_0 r_0^\nu)|=0$ for $\nu\ne l$, $m_l(\chi)m_k(\chi)=0$ when $k\neq l$.\\
 {\bf Step 6.} Define the functions
 $$
  \psi_l(x)=\sum_{h\in H_0^{(N+1)}}\beta_h^{(l)}\varphi({\cal
  A}x\dot-h).
 $$
 \begin{theorem}
 The functions $\psi_l(x\dot-h)$, where $l=\overline{1,p-1}$, $h\in
 H_0$, form an orthonormal basis for $W_0$.
 \end{theorem}
 {\bf Proof.} a) We claim that
 $(\varphi(\cdot\dot-g^{(1)}),\psi_l(\cdot\dot-g^{(2)}))=0$ for any
 $g^{(1)}$,  $g^{(2)}\in H_0$. Since
 $$
  \hat\varphi_{\cdot\dot-h}(\chi)=\overline{(\chi,h)}\hat\varphi(\chi),\;\;\hat\varphi_{{\cal
  A}\cdot\dot-g}(\chi)=\frac{1}{p}\overline{(\chi,{\cal
  A}^{-1}g)}\hat\varphi(\chi{\cal A}^{-1}),
 $$
 it follows that
 $$
  (\varphi(\cdot\dot-g^{(1)}),\psi_l(\cdot\dot-g^{(2)}))=\int\limits_X\hat\varphi(\chi)\overline{\hat\varphi(\chi{\cal
  A}^{-1})}\overline{(\chi,g^{(1)})}(\chi,g^{(2)})\overline{m_l(\chi)}\,d\nu(\chi)=0
 $$
 because ${\rm supp}\,\hat\varphi(\chi)=E$ and $m_l(E)=0$,
 $l=\overline{1,p-1}$.\\
b) By analogy
$$
(\psi_k(\cdot\dot-g^{(1)}),\psi_l(\cdot\dot-g^{(2)}))=\int\limits_X|\hat\varphi(\chi{\cal
A}^{-1})|^2(\chi,g^{(2)}\dot-g^{(1)})m_k(\chi)\overline{m_l(\chi)}\,d\nu(\chi)=0
$$
when $k\ne l$.\\
c) We verify that
$(\psi_l(\cdot\dot-g^{(1)}),\psi_l(\cdot\dot-g^{(2)}))=0$,
provided that $g^{(1)},g^{(2)}\in H_0$ and $g^{(1)}\ne g^{(2)}$.
Write this scalar product in the form
$$
(\psi_l(\cdot\dot-g^{(1)}),\psi_l(\cdot\dot-g^{(2)}))=\int\limits_X|\hat\varphi(\chi{\cal
A}^{-1})|^2(\chi,g^{(2)}\dot-g^{(1)})|m_l(\chi)|^2\,d\nu(\chi)=
$$
$$
=\int\limits_{E{\cal A}\bigcap
X_0r_0^l}(\chi,g^{(2)}\dot-g^{(1)})\,d\nu(\chi).
$$
Show that $E{\cal A}\bigcap X_0r_0^l$ is an $(N,H-2N)$-elementary
set. By the definition
\begin{equation}              \label{eq5.3}
E=\bigsqcup\limits_{\overline{\alpha}\in T(V)}\mathfrak
G_{-N}^\bot r_{-N}^{\alpha_{-N}}...r_{0}^{\alpha_{0}}\dots
r_{s}^{\alpha_{s}}r_{s+1}^{0}...r_{s+N}^{0}\;\;(s\le H-2N-1)
\end{equation}
where the union is taken over all paths
$$
\overline{\alpha}=(0,...,0,\alpha_s,\alpha_{s-1},\dots,\alpha_0,\alpha_{-1},...,\alpha_{-N})
$$
of the tree $T$. It means that for any vector
$(\alpha_{-1},...,\alpha_{-N}),\ \alpha_j=\overline{0,p-1}$ the
union \eqref{eq5.3} contains unique coset $\mathfrak G_{-N}^\bot
r_{-N}^{\alpha_{-N}}...r_{-1}^{\alpha_{-1}}r_{0}^{\alpha_{0}}\dots
r_{s}^{\alpha_{s}}r_{s+1}^{0}...r_{s+N}^{0} $.

Consequently $E{\cal A}=$
$$
 \bigsqcup\limits_{\overline{\alpha}\in T(V)}\mathfrak G_{-N+1}^\bot
r_{-N+1}^{\alpha_{-N}}...r_{0}^{\alpha_{-1}}\dots
r_{s+1}^{\alpha_{s}}r_{s+2}^{0}...r_{s+N+1}^{0}=
$$

$$
\bigsqcup\limits_{\alpha_{-N-1}=0}^{p-1}
\bigsqcup\limits_{\overline{\alpha} \in T(V) }\mathfrak
G_{-N}^\bot r_{-N}^{\alpha_{-N-1}}
r_{-N+1}^{\alpha_{-N}}...r_{0}^{\alpha_{-1}}\dots
r_{s+1}^{\alpha_{s}}r_{s+2}^{0}...r_{s+N+1}^{0}.
$$

On the other hand
$$
X_0r_0^l=\bigcup\limits_{j\in\mathbb
N}\bigsqcup\limits_{(\gamma_{-N},...,\gamma_{-1},\gamma_0)\in
T(V)}\bigsqcup\limits_{b_1,b_2,\dots,b_j=0}^{p-1}\mathfrak
G_{-N}^\bot
r_{-N}^{\gamma_{-N}}...r_{-1}^{\gamma_{-1}}r_{0}^{\gamma_{0}+l}r_1^{b_1}\dots
r_{j}^{b_j}.
$$
Therefore $E{\cal A}\bigcap X_0r_0^l$ consists of all cosets
 $$
  \mathfrak G_{-N}^\bot
   r_{-N}^{\gamma_{-N}}...r_{-1}^{\gamma_{-1}}r_{0}^{\alpha_{-1}}r_1^{\alpha_0}\dots
   r_{s+1}^{\alpha_{s}}r_{s+2}^{0}...r_{s+N+1}^{0}
$$
 where
$$
{(0,...,0,\alpha_s,\alpha_{s-1},\dots,\alpha_{-1}= \gamma_0+l,
\gamma_{-1},...,\gamma_{-N})\in T}
$$
 Since the tree $T$ is $N$-valid it follows that $
 E{\cal A}\bigcap X_0r_0^l$ is $(N,H-2N+1)$-elementary
set. By lemma 3.4 it follows that
$$
\int\limits_{E{\cal A}\bigcap
X_0r_0^l}(\chi,g^{(2)}\dot-g^{(1)})\,d\nu(\chi)=0.
$$
d) We claim that any function $f\in W_0$ can be expanded uniquely
in a series in terms of $(\psi_l(x\dot-g))_{l=\overline{1,p-1},g\in H_0}$.
The proof of this fact may be found in [13], theorem 5.1.
$\square$\\
{\bf Step 7.} Since the subspaces $(V_j)_{j\in\mathbb Z}$ form an
MRA in $L_2(\mathfrak G)$, it follows that the functions
$$
(\psi_l({\cal A}^nx\dot-h))\;\;l=\overline{1,p-1}, n\in\mathbb Z,
h\in H_0
$$
form a complete orthogonal system in $L_2(\mathfrak G)$.


\newpage


 \begin{thebibliography}{99}
  \bibitem{1} Lang W.C., Orthogonal wavelets on the Cantor dyadic group,
    SIAM J.Math. Anal., 1996, 27:1 ,305-312.
  \bibitem{2} Lang W.C., Wavelet
    analysis on the Cantor dyadic group. Housten J.Math.,1998, 24:3, 533-544.
  \bibitem{3} Lang W.C., "Fractal multiwavelets related to
    the Cantor dyadic group, Internat. J. Math. Math. Sci., 1998, 21:2, 307-314.
  \bibitem{4} V Yu Protasov, Y. A. Farkov. Dyadic wavelets and refinable functions on a half-line
 Sbornik: Mathematics(2006), 197(10):1529
  \bibitem{5} Y. A. Farkov, Orthogonalwavelets with compact support on locally compact abelian
    groups, Izvestiya RAN: Ser. Mat., vol. 69, no. 3, pp. 193-220,
    2005, English transl., Izvestiya: Mathematics, 69: 3 (2005), pp.
    623-650.

 \bibitem{6} Y. A. Farkov, Orthogonal wavelets on direct products of
    cyclic groups, Mat. Zametki, vol. 82, no. 6, pp. 934-952, 2007,
    English transl., Math. Notes: 82: 6 (2007).
 \bibitem{7} Yu. Farkov. Multiresolution Analysis and Wavelets
    on Vilenkin Groups. Facta universitatis, Ser.: Elec. Enerd.
    vol. 21, no. 3, December 2008, 309-325
 \bibitem{8} Yu.A. Farkov, E.A. Rodionov. Algorithms for Wavelet Construction on Vilenkin Groups.
    p-Adic Numbers, Ultrametric Analysis and Applications, 2011, Vol. 3, No. 3, pp. 181-195.
  \bibitem{9}  A.\,Yu.~Khrennikov, V.\,M.~Shelkovich, M.~Skopina.
  \it $p$-Adic orthogonal Wavelet Bases.~ \rm P-adic numbers,
  Ultrametric Analysis and Applications, 1:2, 2009,145-156.

  \bibitem{10}  A.\,Yu.~Khrennikov, V.\,M.~Shelkovich, M.~Skopina
  \it $p$-Adic refinable functions and MRA-based wavelets.\rm
  J.Approx.Theory. 161:1, 2009,226-238.

  \bibitem{11} S. Albeverio, S. Evdokimov, M. Skopina
    p-Adic Multiresolution Analysis and Wavelet Frames, J Fourier Anal Appl, (2010), 16: 693-714
  \bibitem{12}  Agaev G.N., Vilenkin N.Ja., Dzafarli G.M., Rubinshtein A.I.,
     Multiplicative systems and harmonic analysis on zero-dimensional groups, ELM, Baku,1981 (in russian).
  \bibitem{13} Lukomskii S.F., Multiresolution analysis on zero-dimensional groups and wavelets bases,
    Math. sbornik, 2010, 201:5 41-64, in russian.
    (english transl.:S.F.Lukomskii, Multiresolution analysis on
    zero-dimensional Abelian groups and wavelets bases,
    SB MATH, 2010, 201:5, 669-691)
  \bibitem{14} Lukomskii S.F. Step refinable functions and orthogonal MRA on $p$-adic Vilenkin
groups. JFAA, February 2014, vol 20, issue 1, pp.42-65.
   \bibitem{15} S. F. Lukomskii. Riesz Multiresolution Analysis on Vilenkin
   Groups. Doklady Mathematics, 2014, Vol. 90, No. 1, pp. 1–4. Original Russian Text
   © S.F. Lukomskii, 2014, published in Doklady Akademii Nauk, 2014, Vol. 457, No. 1,
   pp. 24–27.



\end{thebibliography}
 \end{document}